\input amstex
\documentstyle{amsppt}
\TagsOnRight
\magnification\magstep1
\tolerance=2000

\hoffset1 true pc
\voffset2 true pc
\hsize36 true pc
\vsize50 true pc

\NoBlackBoxes

\define\m1{^{-1}}
\def\gp#1{\langle#1\rangle}

\topmatter
\title
Twisted Group Rings Whose Units Form an $FC$-Group
\endtitle

\author
Victor Bovdi
\endauthor
\thanks
Research supported by the  OTKA Research grant N: T014279
\endthanks
\address
\noindent
Department of Mathematics,
\newline
Bessenyei Teachers' Training College,
\newline
4400, Ny\'\i regyh\'aza,
\newline
 Hungary
\newline
email: vbovdi\@math.klte.hu
\endaddress

\abstract
Let $U(K_\lambda G)$ be the group of units of an infinite twisted group
algebra $K_\lambda G$ over a field $K$. We describe the maximal
$FC$-subgroup of $U(K_\lambda G)$ and give a characterization of
$U(K_\lambda G)$ with finitely conjugacy classes. In the case of
group algebras we obtain the Cliff-Sehgal-Zassenhaus' theorem.
\endabstract
\endtopmatter

\document
\subhead 1. Introduction \endsubhead Let $G$ be a group, $K$ a
field and $\lambda:G\times G\longmapsto U(K)$ a $2$-cocycle of
$G$ with respect to the trivial action $G$. Then the twisted group
algebra $K_{\lambda}G$ of $G$ over field $K$ is an associative
$K$-algebra with basis $\{u_g \mid g\in G\}$ and with
multiplication defined for all $g,h\in G$
$$
u_gu_h=\lambda_{g,h}u_{gh},\qquad (\lambda_{g,h}\in \lambda)
$$
and using distributivity.

Let $U(K_\lambda G)$ be the group of units of $K_\lambda G$ and  let
$\Delta U$ be its  subgroup  consisting of all elements with  finitely
many conjugates in $U(K_\lambda G)$.  This subgroup $\Delta U$ is called
the $FC$-center of $U(K_\lambda G)$.  Clearly, if $\Delta
U=U(K_\lambda G)$, then  $U(K_\lambda G)$ is an $FC$-group (group
with finite conjugacy classes).

The problem to study the group of units of group rings  with $FC$
property was posed by S.K.~Sehgal and H.J.~Zassenhaus \cite{1}.
For a field $K$ of characteristic $0$ they described  all groups  $G$ without
subgroups of type $p^\infty$ for which the group of units of the group algebra
of $G$ over $K$ is an $FC$ group. This was spelling for arbitrary groups by H.Cliff
and S.K.Sehgal \cite{2}.

In this paper we describe the subgroup  $\Delta U$ when $K_\lambda G$ 
is infinite. Let $t(\Delta U)$ be the group of  all elements of finite order of $\Delta
U$. Then $\Delta U$ is a solvable group of length at most $3$ and
subgroup $t(\Delta U)$ is of nilpotent class at most $2$. This is new even for group algebras.
We use this result for characterization of those cases when $U(K_\lambda G)$ has $FC$ property, and
obtain a generation  of the Cliff-Sehgal-Zassenhaus theorem for
twisted group algebras.

\subhead 2. The $FC$ center of $U(K_\lambda G)$ \endsubhead By the
Theorem of B.H. Neumann \cite{3} the elements of finite order in
$\Delta U$ form a normal subgroup which we denote by $t(\Delta
U)$, and the factorgroup $\Delta U / t(\Delta U)$ is a torsion free
abelian group.  Evidently, $\overline{G}=\{\kappa u_a\mid \kappa
\in U(K), a\in G\}$ is a subgroup in $U(K_\lambda G)$, while
$U(K)$ is a normal subgroup in $\overline{G}$, with factorgroup
$\overline{G}/U(K)$ isomorphic to $G$.

If $x$ is a nilpotent element of the ring $K_\lambda G$ then the
element $y=1+x$ is a unit in $K_\lambda G$ and is referred to as a
unipotent element of $U(K_\lambda G)$.

Let $\zeta (G)$ be the center of group $G$ and $[g, h]=g^{-1}h^{-1}gh$
\qquad ($g, h\in G$).

\proclaim{Lemma 1}
{\it Let $K_\lambda G$ be an infinite twisted group
algebra. Then all unipotent elements of the subgroup $\Delta U$ are
central in $\Delta U$.}
\endproclaim
{\bf Proof.}
Let $y=1+x$ be a unipotent element of $\Delta U$ and $v\in
\Delta U$. Then for a positive integer $k$ we have $x^k=0$ and
by induction on $k$ we will prove $vx=xv$.

The subgroup $\overline{G}=\{\kappa u_a\mid \kappa
\in U(K), a\in G\}$ is infinite and by Poincar\'e's Theorem
the centralizer $S$ of the subset $\{v, y\}$ of $\overline{G}$ is
a subgroup of finite index in $\overline{G}$. Since
$\overline{G}$ is infinite, $S$ is infinite and $fy=yf$ for all
$f\in S$. Then $xf$ is nilpotent and $1+xf$ is a unit in
$K_\lambda G$. We can see easily that the set
$\{(1+xf)^{-1}v(1+xf) \mid f\in S \}$ is finite. Let $v_1, \ldots,
v_t$ be the elements of this set and $W_i=\{f\in S\mid
(1+xf)^{-1}v(1+xf)=v_i\}$. Then $S=\cup W_i$ and there exists an
index $j$ such that $W_j$ is infinite. Fix an element $f\in W_j$.
Any element $q\in W_j$, $q\not=f$ satisfies
$$
(1+xf)^{-1}v(1+xf)=(1+xq)^{-1}v(1+xq)
$$
and
$$
v(1+xf)(1+xq)^{-1}=(1+xf)(1+xq)^{-1}v.
$$
Then
$$
\split
v\{(1+xq)+(xf-xq)\}(1+xq)^{-1}&=\{(1+xq)+(xf-xq)\}(1+xq)^{-1}v,\\
v(1+x(f-q)(1+xq)^{-1})&=(1+x(f-q)(1+xq)^{-1})v
\endsplit
$$
and
$$
vx(f-q)(1+xq)^{-1}=x(f-q)(1+xq)^{-1}v.\tag1
$$
Now we use the induction mentioned above.  For $k=1$ the statement is trivial; so we suppose it is true for all
$1\leq n<k$, where $k\geq 2$ is any given integers.

If $m\geq 2$, then by induction hypothesis $x^mv= vx^m$ for all $v\in \Delta U$.
Clearly, if $i\geq 1$ then
$$
x(f-g)x^{i}q^iv=(f-g)x^{i+1}q^iv=(f-g)vx^{i+1}q^i=vx(f-g)x^{i}q^i.
$$
From (1) we have
$$
\split
vx(f-q)(1-xq&+x^2q^2+\cdots+(-1)^{k-1}x^{k-1}q^{k-1})\\
&=x(f-q)(1-xq+x^2q^2+\cdots+(-1)^{k-1}x^{k-1}q^{k-1})v.
\endsplit
$$
So $(f-q)(vx-xv)=0$.

Now suppose $xv\not=vx$. The element
$q^{-1}f \in\overline{G}$ can be written as $\lambda u_h$
\newline
($\lambda\in U(K), h\in G$). By $vx-xv=\sum_{i=1}^{s}
\alpha_iu_{g_i}\not=0$, then we have
$$
\sum_{i=1}^{s}\lambda\alpha_iu_hu_{g_i}-\sum_{i=1}^{s}\alpha_iu_{g_i}=0.
$$
If $h\in G$ satisfies this equation, then $g_1=hg_j$ for some $j$, and the
number of such elements $h$ is finite. Since $W_j=\{\lambda u_h \mid
\lambda \in U(K) \}$ is an infinite set, there exist $h$ and
different elements $\lambda_1, \lambda_2\in K$ such that
$\lambda_1u_h, \lambda_2u_h\in W_j$.
Then $(\lambda_i u_h-1)(vx-xv)=0$, \quad ($i=1, 2$) and we obtain
$(\lambda_1u_h-\lambda_2u_h)(vx-xv)=0$. This condition is satisfied  only
if $vx=xv$ but does not hold. \hfil\qed

\proclaim{Lemma 2}
{\it Let $K_\lambda G$ be an infinite twisted group algebra, $H$  a
finite subgroup of $\Delta U$ and $L_H$ the subalgebra of
$K_\lambda G$ generated by $H$. Then the group of units $U(L_H)$ of
the algebra $L_H$ is contained in $\Delta U$, and the factorgroup
$U(L_H)/(1+J(L_H))$ is abelian.}

\endproclaim
{\bf Proof.} If $H$ is a finite subgroup of $\Delta U$ and $L_H$ is the subalgebra
of $K_\lambda G$ generated by $H$, then $L_H$ is a finite algebra
over $K$ and its radical $J(L_H)$ is nilpotent. Then $U(L_H)$ is a
subgroup of $\Delta U$ and by Lemma $1$ all unipotent elements of
$U(L_H)$ are central in $\Delta U$. Therefore $1+J(L_H)$ is a central
subgroup of $\Delta U$ and $J(L_H)\subset\zeta (L_H)$, where $\zeta
(L_H)$ is the center of $L_H$. Then by Theorem 48.3 in \cite{4} (p.209)
$$
L_H=L_He_1\oplus \cdots \oplus L_He_n\oplus N, \tag2
$$
where $L_He_i$ is a semiprime algebra (i.e. $L_He_i/J(L_He_i)$ is
a division ring), $N$ is a commutative artinian radical algebra,
$e_1, \ldots,  e_n$ are pairwise orthogonal idempotents. By Lemma
$13.2$ in \cite{4} $(p.57)$ any idempotent $e_i$ is central in
$L_H$ and $U(L_He_i)$ is isomorphic to the subgroup
$\gp{1-e_i+ze_i \mid z\in U(L_H)}$ of $U(L_H)$.

Since $U(L_He_i)$ is a subgroup of the $FC$-group $\Delta U$, it is an $FC$-group, too.
As $J(L_He_i)$ is nilpotent (see \cite{5}),
$$
                U(L_He_i)/(1+J(L_He_i))\cong U(L_He_i/J(L_He_i)). \tag3
$$
By Scott's Theorem \cite{7}, in the skewfield $L_He_i/J(L_He_i)$ every nonzero element
is either central or its conjugacy class is infinite. Thus the $FC$-group
$U(L_He_i)/(1+J(L_He_i))$ is abelian.

Decomposition ($2$) implies
$$
         L_H/J(L_H)\cong L_He_1/J(L_He_1)\oplus\cdots \oplus L_He_n/J(L_He_n)
$$
and
$$
\split U(L_H)/(1+J(L_H))&\cong U(L_H/J(L_H))\\
& \cong U(L_He_1/J(L_He_1))\times \cdots \times U(L_He_n/
J(L_He_n)).
\endsplit
$$
Therefore $U(L_H)/(1+J(L_H))$ is abelian. \hfill \qed

\proclaim {Theorem 1}
{\it Let $K_\lambda G$ be an infinite twisted group algebra and
$t(\Delta U)$  the subgroup of $\Delta U$ consisting of  all elements
of finite order in $\Delta U$. Then all elements of the commutator
subgroup of $t(\Delta U)$ are unipotent and central in $\Delta
U$.}
\endproclaim

{\bf Proof.}
Let $H$ be a finite subgroup of $t(\Delta U)$ and $L_H$ be the
subalgebra of $K_\lambda G$, generated by $H$. Then the elements of
the subgroup $H_1=H\cap (1+J(L_H))$ are unipotent and (by Lemma 1)
central in $\Delta U$. The subgroup $H(1+J(L_H))$ is contained in
$U(L_H)$ and
$$
H/H_1=H/(H\cap (1+J(L_H)))\cong
(H(1+J(L_H)))/(1+J(L_H)).
$$
By Lemma 2 the factorgroup
$U(L_H)/(1+J(L_H))$ is abelian. So $H/H_1$ is abelian and the
commutator subgroup of $H$ is contained in $H_1$ and consists of
unipotent elements.

Since the commutator subgroup of $t(\Delta U)$ is the union of the
commutator subgroups of the finite subgroups of $t(\Delta U)$, all
elements of the commutator subgroup of $t(\Delta U)$ are unipotent
and are central in $\Delta U$.
\hfill \qed

\proclaim{Theorem 2}
{\it Let $K_\lambda G$ be an infinite twisted group algebra such
that $char(K)$ does not divide the order of any element of the
subgroup $\Delta G$. Then $t(\Delta U)$ is abelian.}
\endproclaim
{\bf Proof.} Let $H$ be a finite subgroup of the commutator
subgroup of $t(\Delta U)$. Then (by Theorem $1$) $H$ is
contained in the center of $\Delta U$. The set $\{
u_g^{-1}Hu_g\mid g\in \Delta G \}$ contains only a finite number
of subgroups $H_1, H_2, \ldots,  H_t$. The subgroup $L=H_1\cdot
H_2\cdots H_t$ is finite and is invariant under the inner automorphism
$f_g(x)=u_g^{-1}xu_g$ of the ring $K_\lambda \Delta G$, where $g\in \Delta G$.
Let $x_1,\ldots, x_s$ be all elements of $L$. Then $y_i=x_i-1$ is a nilpotent
element, and in the commutative ring $L$ the elements $y_1,
\ldots,  y_s$ commute.  Therefore
$$
J\cong\{\sum_{i=1}^{s} \alpha_iy_i \mid \alpha_i\in K, x_i=y_i+1\in L\}
$$
is a nilpotent subring. Let
$$
F=\{\sum_{i=1}^{s}\alpha_iy_iz_i
\mid\alpha_i\in K, x_i=y_i+1\in L, z_i\in K_\lambda \Delta G\}.
$$
Let us prove that $F$ is a nilpotent right ideal of $K_\lambda
\Delta G$. If $z=\sum_{j} \beta_ju_{g_j}\in K_\lambda \Delta G$
then $y_iz=\sum_{j}\beta_ju_{g_j}u_{g_j}^{-1}y_iu_{g_j}$, and
$u_{g_j}^{-1}y_iu_{g_j}$ equals one of the elements $ y_1, \ldots,
y_s$. This and the nilpotency of the ring $J$ imply that $F$ is a
nilpotent ring. By Passman's Theorem \cite{6}, if $char(K)$ does
not divide the order of any element of $\Delta G$ then $K_\lambda
\Delta G$ does not contain nilideals. Therefore $F=0$, $L=1$ and
the commutator subgroup $t(\Delta U)$ is trivial, so $t(\Delta U)$
is abelian. \hfill \qed

\proclaim{Corollary}
{\it Let $K_\lambda \Delta G$ be an infinite twisted group algebra. Then
$\Delta U$ is a solvable group of length at most $3$, and subgroup
$t(\Delta U)$ is of nilpotent class at most $2$}.
\endproclaim

\subhead 3. The $FC$ property of $U(K_\lambda G)$ \endsubhead

\proclaim{Lemma 3} { \it Let $L$ be a subfield of the twisted group
algebra $K_\lambda G$,  where $K$ is a subfield of $L$,  $g\in G$
an element of order $n$ and
$$
\lambda_g=u_g^n=\lambda_{g, g}\lambda_{g, g^2}\cdots \lambda_{g, g^{n-1}}.
$$
If $\alpha^n \not=\lambda _g$ for some $\alpha \in L$ and
$\alpha u_g=u_g\alpha$ then $u_g-\alpha$ is a unit in $K_\lambda G$.
Furthermore, if $L$ is an infinite field then the number of such
units is infinite.

}
\endproclaim
{\bf Proof.} Let $\alpha \in L, \alpha^n\not=\lambda _g$ and
$u_g\alpha=\alpha u_g$. Then $\lambda _g-\alpha^n$ is a nonzero
element of L and
$$
(\alpha^{n-1}+\alpha^{n-2}u_g+\cdots +\alpha
u_g^{n-2}+u_g^{n-2})(\lambda _g-\alpha ^n)^{-1}
$$
is the inverse of $u_g-\alpha$. We know that the number of solutions
of the equation $x^n-\lambda_g=0$ in L does not exceed $n$. Thus in
an infinite field $L$ there are infinitely many elements not
satisfying the equation $x^n-\lambda_g=0$.
\hfill \qed

\proclaim{Lemma 4}
{ \it Let $G$ be an infinite locally finite group and  $char(K)$ does
not divide the order of any element of $G$. If $U(K_\lambda G)$ is an
$FC$-group then $G$ is abelian and $K_\lambda G$ is commutative.

}
\endproclaim
{\bf Proof.} Let $W$ be a finite subgroup of $G$. Then the
subalgebra $K_\lambda W$ is a semiprime artinian ring and by the
Wedderburn-Artin Theorem
$$
K_\lambda W=M(n_1, D_1)\oplus \cdots \oplus M(n_t, D_t),
$$
where each $D_k$ is a skewfield and $M(n_k, D_k)$ is a full matrix algebra.
Let $e_{i, j}, e_{j, i}$ be  matrix units in $M(n_k,
D_k)$ and $i\not=j$. Then the unipotent elements $1+e_{i, j}$,
$1+e_{j, i}$ are central in $K_\lambda G$ (see Theorem 1) which is
impossible if $i\not=j$. Thus $n_k=1$ and $K_\lambda W$ is a
direct sum of skewfields, $K_\lambda W=D_1\oplus D_2\oplus \cdots
\oplus D_t$ and
$$
          U(K_\lambda W)=U(D_1)\times U(D_2)\times \cdots \times U(D_t).
$$
By Scott's Theorem \cite{7} any nonzero element of a skewfield is either
central or has an infinite number of conjugates. Therefore $K_\lambda W$
is a direct sum of fields and $W$ is abelian. Since $G$ is a
locally finite group, $G$ is abelian and $K_\lambda
G$ is a commutative algebra.
\hfill \qed

\proclaim{Lemma 5}
{\it Let $K_\lambda G$ be infinite and $char(K)$ does not divide the order of any element of the normal torsion subgroup  $L$ of $G$. If $U(K_\lambda G)$ is an $FC$-group, then all idempotents of $K_\lambda L$ are cental
in $K_\lambda L$.}
\endproclaim
{\bf Proof.} Let the idempotent $e\in K_\lambda L$ be noncentral in $K_\lambda G$.
Then there are exists $g\in G$ such that $eu_g\not=u_ge$. The subgroup
$H=\gp{g^{-i}supp(e)g^i\mid i\in \Bbb Z}$ is finite and for any $a\in G$ the subalgebra
$K_\lambda H$ of $K_\lambda L$ is invariant under the inner automorphism $\phi(x)=u_a\m1xu_a$. It is easy to see (by Lemma 4) that $K_\lambda H$ is a commutative semisimple $K$-algebra of finite rank and the idempotent $e\in K_\lambda H$ is a sum of primitive idempotents. Consequently, there exists a primitive idempotent $f$ of $K_\lambda H$ which does not commute with $u_g$.  Then the idempotents $f$ and $u_g\m1fu_g$ are orthogonal  and $(u_gf)^2=u_gfu_gf=u_g^2(u_g\m1fu_g)f=0$. By Theorem 1 the unipotent element $1+u_gf$ commutes with $u_g$ and
$(1+u_gf)u_g=u_g(1+u_gf)$ implies $u_gf=fu_g$, which is impossible. Thus, all idempotents of $K_\lambda L$ are central in $K_\lambda G$.\hfill \qed

\proclaim{Lemma 6}
{ \it Let $U(K_\lambda G)$ be an $FC$-group and $t(G)$ the set of
elements of finite order in $G$. Then
\item{1.} $G$ is an $FC$-group;
\item{2.} if there exists an infinite subfield $L$ in  the center of $K_\lambda G$ such
that $L\supseteq K$ then $t(G)$ is central in G and $\lambda_{g,
h}=\lambda_{h, g}$ \quad  ($h\in t(G)$, $g\in G$).

}
\endproclaim
{\bf Proof.}
If $U(K_\lambda G)$ is an $FC$-group then $\overline{G}=\{
\lambda u_g\mid\lambda \in U(K), g\in G\}$ is an $FC$-subgroup. Clearly
$U(K)$ is normal in $\overline{G}$ and $G\cong\overline{G}/U(K)$. We conclude
that $G$ is an $FC$-group as it is a homomorphic image of the $FC$-group
$\overline{G}$.

Let $L$ be an infinite field which satisfies condition 2. of the
Lemma. Then by Lemma 1 for any $h\in t(G)$ there exists a
countable set $S=\{\alpha_i \in L\mid i\in \Bbb Z\}$ such that
$u_h-\alpha_i$ is a unit for all $i\in \Bbb Z$. Suppose that
$u_gu_h\not=u_hu_g$ for some $g\in G$. Next we observe that the
equality
$$
                   (u_h-\alpha_i)u_g(u_h-\alpha_i)^{-1}=
                    (u_h-\alpha_j)u_g(u_g-\alpha_j)^{-1}
$$
holds  only in case $\alpha_i=\alpha_j$. Since
$$
                        (u_h-\alpha_i)(u_h-\alpha_j)^{-1}=
                     1+(\alpha_j-\alpha_i)(u_h-\alpha_j)^{-1},
$$
we obtain that $(\alpha_i-\alpha_j)(u_gu_h-u_hu_g)=0$ and
$\alpha_i=\alpha_j$. It follows that the set

$$
        \{(u_h-\alpha_j)u_g(u_h-\alpha_j)^{-1} \mid i\in \Bbb Z \}
$$
is infinite which contradicts the condition that $U(K_\lambda G)$ is an
$FC$-group. Then $u_gu_h=u_hu_g$, so $[g, h]=1$, $t(G)\subseteq
\zeta(G)$ and $\lambda_{g, h}=\lambda_{h, g}$\quad  ($h\in t(G), g\in G$).
\hfill \qed

\proclaim{Lemma 7}
{\it Let G be an  abelian torsion group, $K_\lambda
G$ a commutative semisimple algebra and $v$ an idempotent of
$K_\lambda G$. If $K_\lambda Gv$ contains a finite number of
idempotents then $K_\lambda Gv$ is a direct sum of finitely many
fields.}
\endproclaim
{\bf Proof.} If $e_1, \ldots,  e_s$ are all the idempotents of
$K_\lambda Gv$, then
$$
L=\gp{supp(e_1), \ldots,  supp(e_t)}
$$
is a finite subgroup in G and $K_\lambda Lv$ is a direct sum of
finitely many fields,
$$
K_\lambda Lv=(K_\lambda Lv)f_1\oplus \cdots
\oplus (K_\lambda Lv)f_t,
$$
where $f_1, \ldots,  f_t$ are orthogonal primitive idempotents of
$K_\lambda Lv$.  The corresponding direct sum in $K_\lambda Gv$ is
$$
K_\lambda Gv=(K_\lambda Gv)f_1\oplus
\cdots \oplus (K_\lambda Gv)f_t.
$$
We show that the element $0\not= x\in (K_\lambda Gv)f_i$ is a
unit. $R=\gp{L, supp(x)}$ is a finite subgroup and $K_\lambda Rv$
is a direct sum of finitely many fields,
$$
        K_\lambda Rv=(K_\lambda Rv)l_1\oplus \cdots \oplus (K_\lambda Rv)l_t,
$$
and each idempotent $f_i$ is either equal to an idempotent $l_j$ or is a sum
of these idempotents. If $f_i=l_j$ then $xf_i\in (K_\lambda Rv)l_j$ and $x$
is a unit in $(K_\lambda Lv)f_i$. If $f_i=l_{i_1}+l_{i_2}$\quad
$(l_{i_1}, l_{i_2}\in K_\lambda Lv)$\quad  then $(K_\lambda Lv)f_i=(K_\lambda
Lv)l_{i_1}\oplus (K_\lambda Lv)l_{i_2}$, but this  does not hold.
\hfill \qed

\proclaim{Theorem 3}
{\it Let $K_\lambda G$ be an infinite twisted group
algebra of $char(K_\lambda G)=p$, such that $t(G)$ contains a $p$-element and
either the field $K$ is perfect or for any element $g\in G$ of order $p^k$,
the element $u_g^{p^k}$ is an algebraic element over the prime subfield of
$K$. Then $U(K_\lambda G)$ is an $FC$-group if and only if $G$ is an
$FC$-group and satisfies the following conditions:

\item{1.} $p=2$ and $\mid G' \mid =2$;

\item{2.} $t(G)$ is central in G and $t(G)=G'\times H$, where $H$
is abelian, and has no $2$-elements;

\item{3.} $K_\lambda H$ is a direct sum of a finite number of fields;
\footnote {If  $K_\lambda H$ is a group ring then $H$ is a finite
abelian group.}

\item{4.} $\{\lambda_{h, h^{-1}}^{-1} \lambda_{h^{-1}, g}\lambda_{h^{-1}g, h}
\mid h\in H \}$ is a finite set for all $g\in G$.}
\endproclaim
{\bf Proof. Necessity.}
By Lemma 6 $G$ is an $FC$-group. Let $g$ be an element of order $p^k$. Then $u_g^{p^k}=\lambda_g\in
U(K)$, and in the perfect field $K$ we can take the $p^k$-th root of
$\lambda_g$ which we denote by $\mu$. If $K_0$ is the prime subfield of $K$ and $\lambda_g$
is algebraic over $K_0$ then $K_0(\lambda_g)$ is a finite field
and so it is perfect. Thus $u_g-\mu$ is nilpotent and
$1+\mu-u_g$ and (by Theorem 1) the element $1-(u_g-\mu )u_a$ are central in $U(K_\lambda G)$.
Then for any $b\in G$ by
$$
                   u_b(1-(u_g-\mu )u_a)=(1-(u_g-\mu )u_a)u_b
$$
implies
$$
                u_bu_gu_a-\mu u_bu_a-u_gu_au_b+\mu u_au_b=0.\tag4
$$
Each  $u_g$ can be written in the form $\mu +(u_g-\mu )$ and so
$\mu^{-1}u_g=1+\mu^{-1}(u_g-\mu )$. Thus $\mu^{-1}u_g$ is an
unipotent element and it commutes with $u_b$ and $u_a$. Then (4)
can be written as
$$
u_gu_bu_a-u_gu_au_b-\mu u_bu_a+\mu u_au_b=0. \tag5
$$
If $[a, b]=1$ then by (5), $(\lambda_{a, b}-\lambda_{b,
a})(u_g-\mu )=0$. From this equation we get that the coefficient
of $u_g$ must be zero and $\lambda_{a, b}=\lambda_{b, a}$. Thus,
$u_bu_a=u_au_b$.

Let $[a, b]\not=1$. Then by (5), $u_gu_bu_a=-\mu u_au_b$ and
$u_gu_au_b=-\mu u_bu_a$. So
$$
\cases
      \text{$u_g=-\mu [u_a^{-1}, u_b^{-1}]^{-1}$,}\\
      \text{$u_g=-\mu[u_a^{-1}, u_b^{-1}]$.}
\endcases                                                            \tag6
$$
Consequently $u_g^2=\mu^2$ and $(u_g\mu^{-1})^2=1$. Note that in
($6$) $g$ may be any $p$-element, further $a$ and $b$ may be any
noncommuting elements of $G$. This is possible only if $p=2$. Then
the commutator subgroup $\overline{G}'$ of group $\overline{G}$ is
of order $2$ and coincides with the Sylow $2$-subgroup of
$\overline{G}$. Thus $\overline{G}'\subseteq \zeta (\overline{G})$
and $\overline{G}$ is a nilpotent group of class at most $2$. Let
$$
             L= \gp{\mu u_h\mid\mu \in U(K), h\in t(G)}.
$$
Then $L/ U(K)$ is a torsion nilpotent group and its $2$-Sylow
subgroup is of order $2$. Here $L$ is abelian, because $\overline{G}'$ is of order $2$
and it is a subgroup in $L$. Therefore $t(G)$
is abelian and $t(G)=S\times H$, where $S=\gp{g\mid g^2=1}$ is the
Sylow $2$-subgroup of $t(G)$ and all elements of $H$ are of odd order.

We show that $K_\lambda H$ is central in $K_\lambda G$. Let $h\in
H, a\in G$ and $[u_a, u_h]\not=1$. Then $[u_a, u_h]=\mu u_g$ and
$$
             \lambda u_{a^{-1}h^{-1}ah}=\mu u_g.                     \tag7
$$
It is clear that  $[a, h]\in H$ and the
order of $[a, h]$ is odd, because $H$ is normal in $G$. Since $g$ is a $2$-element,
(7) does not hold. Thus $K_\lambda H$ is central in $K_\lambda
G$ and $t(G)\subseteq \zeta (G)$.

Let us prove that $K_\lambda H$ contains only a finite number of
idempotents. Suppose $K_\lambda H$ contains an infinite number of
idempotents $e_1, e_2, \ldots $. If $d, b\in G$ and $[b,
d]=g\not=1$ then $g^2=1$ and (by Lemma 5) the element $1-e_i+u_de_i$ is a unit. Clearly,
$$
         (1-e_i+u_de_i)^{-1}u_b(1-e_i+u_de_i)=u_b(1-e_i+\mu u_ge_i),
$$
where $\mu =\lambda_{d, d^{-1}}^{-1}\lambda_{b, b^{-1}}^{-1}
\lambda_{d^{-1}, b} \lambda_{d^{-1}b, d} \lambda_{d^{-1}bd, b^{-1}}$.

If $i\not=j$ then $1-e_i+\mu u_ge_i\not=1-e_j+\mu u_ge_j$. Indeed,
if
$$
1-e_i+\mu u_ge_i=1-e_j+\mu u_ge_j,
$$
then $(e_i-e_j)(\mu
u_g-1)=0$. Since $e_i-e_j\in K_\lambda H$ and $u_g\notin K_\lambda
H$, the last equality is true only in case $i=j$. Therefore if
$i\not=j$ then $1-e_i+\mu u_ge_i\not=1-e_j+\mu u_ge_j$ and $u_b$
has an infinite number of conjugates, which does not hold. Thus
$K_\lambda H$ contains a finite number of idempotents $e_1,
\ldots,  e_t$, and (by Lemma 7) $K_\lambda H$ is a direct sum of a
finite number of fields.

Since $\{u_g^{-1}u_hu_g \mid g\in G\}$ is a finite set, we obtain condition 4.
of the Theorem.

\smallskip
{\bf Sufficiency.}
Let the conditions of the Theorem be satisfied.  We prove that
$U(K_\lambda G)$ is an $FC$-group.

Let $G'=\gp{a\mid a^2=1}$ be the commutator subgroup of $G$ and
$\mu^2=\lambda_{a, a}$. Thus the ideal $\frak I=K_\lambda G(u_a-\mu)$
is nilpotent.

In $K_\lambda G$ we choose a new basis $\{w_g\mid g\in G\}$,
$$
w_g=\cases
              u_g,          &\text{if   }  g\in G\setminus \gp{a}, \\
              \mu^{-1}u_g,  &\text{if   }  g\in \gp{a}.
    \endcases
$$
Let $G=\cup b_j\gp{a}$ be the decomposition of the group G by the
left cosets of $\gp{a}$. The element $x+\frak I\in K_\lambda
G/\frak I$ can be written as
$$
\split x+\frak I&=\sum_{i} \alpha_iw_{b_i}+\sum_{i}
\beta_iw_{b_i}w_a+\frak I\\
&= \sum_{i}\alpha_i w_{b_i}+\sum_{i} \beta_iw_b(w_a-1)
+\sum_{i}\beta_iw_{b_i}+\frak I=\sum_{i}( \alpha_i+ \beta_i)w_{b_i}+\frak I.
\endsplit
$$
We show that $K_\lambda G/\frak I$ is commutative. Indeed
$$
(w_g+\frak I)(w_h+\frak I)=\- w_gw_h+\frak I=\-w_hw_g[w_g, w_h]+\frak I,
$$
and the commutator $[w_g, w_h]$ is either 1 or $w_a$. If $[w_g, w_h]=w_a$ then
$$
w_gw_h+\frak I=w_hw_gw_a+\frak I=w_hw_g(w_a-1)+w_hw_g+\frak I=w_hw_g+\frak I.
$$
We will construct the twisted group algebra $K_\mu H$ of the group $H=G/\gp{a}$
over the field $K$ with the system of factors $\mu$.

Let $R_l(G/\gp{a})$ be a fixed set of representatives of all left
cosets of the subgroup $\gp{a}$ in $G$ and $H=\gp{h_i=b_i\gp{a}
\mid b_i\in R_l(G/\gp{a})}$. The element $w_{b_i}+\frak I$ is
denoted by $t_{h_i}$. If $h_ih_j=h_k$, then $b_ib_j=b_ka^s$
($s=\{0, 1\}$), and
$$
\split t_{h_i}t_{h_j}&= w_{b_i}w_{b_j}+\frak I= \lambda_{b_i,
b_j}w_{b_ka^s}+\frak I=\lambda_{b_i, b_j}\lambda_{b_k, a^s}^{-1}
w_{b_k}w_{a^s}+\frak I\\
& =\lambda_{b_i, b_j} \lambda_{b_k, a^s}^{-1}
w_{b_k}+\lambda_{b_i, b_j} \lambda_{b_k, a^s}^{-1}
w_{b_k}(w_{a^s}-1)+\frak I
\\
&=\lambda_{b_i, b_j} \lambda_{b_k, a^s}^{-1} w_{b_k}+\frak I.
\endsplit
$$
Let $\mu_{h_i, h_j}=\lambda_{b_i, b_j}\lambda_{b_k, a^s}^{-1}$ and $\mu
=\{\mu_{a, b} \mid a, b \in H\}$. Let $\{t_h \mid h\in H\}$ be a basis of
the twisted group algebra $K\mu H$ with the system of factors $\mu$.
Clearly,  $t_{h_i}t_{h_j}=\mu_{b_i, b_j}t_{h_k}$.

Let $t(H)$ be the set of elements of finite order of $H$ and
$H=\cup c_it(H)$ the decomposition of the group $H$ by the cosets
of the subgroup $t(H)$. Then elements $x, x^{-1}\in U(K_\mu H)$
can be given as
$$
               x=\sum_{i=1}^{t} \alpha_it_{c_i} \quad
            \text{and} \quad  x^{-1}=\sum_{i=1}^{s}\beta_it_{d_i},
$$
where $\alpha_i, \beta_j$ are nonzero elements of $K_\mu t(H)$.
The subgroup
$$
   L=\gp{supp(\alpha_1), \ldots,  supp(\alpha_t),
   supp(\beta_1), \ldots,  supp(\beta_s)}
$$
is finite and $K_\mu L$ is a direct sum of fields

$$
            K_\mu L=e_1K_\mu L\oplus \cdots \oplus e_nK_\mu L. \tag8
$$
Let $xe_k=\sum_{i=1}^{n} \gamma_it_c{_i}$ and $x^{-1}e_k=\sum_{i=1}^{m}
\delta_it_{d_i}$, where $\gamma_i, \delta_j$ are nonzero elements
of the field $K_\mu Le_k$.

We know \cite{8}, that a torsionfree abelian  group is orderable. Therefore we
can assume that
$$
                    c_{i_1}t(H)<c_{i_2}t(H)<\cdots <c_{i_n}t(H)
$$
and
$$
                    d_{j_1}t(H)<d_{j_2}t(H)<\cdots <d_{j_m}t(H).
$$
Then $c_{i_1}d_{j_1}t(H)$ is called the least and $c_{i_n}d_{j_m}t(H)$ is
called the greatest among the elements of the form $c_{i_s}d_{j_q}t(H)$. It
is easy to see that $c_{i_1}d_{j_1}t(H)<c_{i_n}d_{j_m}t(H)$ if $n>1$ or
$m>1$.  Therefore $\gamma\delta_1t_{c_{i_1}}t_{d_{j_1}}
\not=\gamma_n\delta_m t_{c_{i_n}}t_{d_{j_m}}$.  Since $x^{-1}e_kxe_k=e_k$,
we have $n=m=1$, $xe_k=\gamma t_{c_r}$ and $x^{-1}e_k= \gamma^{-1}
t_{c_r}^{-1}$. Thus, elements $x$ and $x^{-1}$ can be written  as
$$
     x=\sum_{i=1}^{t}\gamma_it_{c_i} \quad \text{ and} \quad
           x^{-1}=\sum_{i=1}^{t}\gamma_i^{-1}t_{c_i}^{-1},
$$
where $\gamma_1, \ldots, \gamma_t$ are orthogonal elements.

Let $\phi :K_\lambda G/\frak I \mapsto K_\mu H$ be an isomorphism
of these algebras. If $x\in U(K_\lambda G)$ then $\phi (x+\frak
I)= \sum_{i=1}^{t} \gamma_it_{c_i}$ and $\gamma_i\in K_\mu Le_i$.
It is easy to see that there exists an abelian subgroup
$\overline{L}$ of $G$ such that $L=\overline{L}/\gp{a}$. The
algebra $K_\lambda \overline{L}$ is commutative and its radical is
a nilpotent ideal equal to $\frak I\cap K_\lambda \overline{L}$.
Since $K_\mu \overline{L}/(\frak I \cap
K_\lambda\overline{L})\cong K_\lambda L$, by the classic method of
lifting idempotents, there exist idempotents $f_1, \ldots,  f_t$
in $K_\mu\overline{L}$ such that $f_1+\cdots +f_t=1$ and
$f_i+\frak I=e_i$. Then $x=xf_1+\cdots +xf_t$ and $\phi
(xf_i+\frak I)=\gamma_it_{c_i}$, where $h_i=b_i\gp{a}, b_i\in G$.
Then there exists an element $v_i\in K_\lambda \overline{L}f_i$
such that $\phi (v_i+\frak I)=\gamma_i$ and $\phi(v_iw_{g_i}+
\frak I)= \gamma_it_{h_i}$. We can find an element $r\in \frak I$
such that $xf_i=(v_i+rf_i)w_{g_i}$.

Clearly  $s_i=v_i+rf_i$ is a unit in $K_\mu \overline{L}f_i$
and is central in $K_\lambda G$. Then $s_1, \ldots,  s_t$ are
orthogonal and \quad $x=\sum_{i=1}^{t}s_iw_{g_i}$,\quad
$x^{-1}=\sum_{i=1}^{t}s_i^{-1}w_{g_i}^{-1}$. Since $s_i\in \zeta
(K_\lambda G)$,\quad  $x^{-1}w_gx=\sum_{i=1}^{t}w_{g_i}^{-1}w_g w_{g_i}$
for any $g\in G$. Because $G$ is an $FC$-group, then by condition
4. of the Theorem $w_g$ has finite number of conjugates. Thus
$U(K_\lambda G)$ is an $FC$-group. \hfill \qed

\proclaim{Lemma 8}
{\it Let $K$ be a field such that $char(K)$ does
not divide the order of any element of $t(G)$, $K_\lambda t(G)$
a commutative algebra that does not contain a nontrivial minimal
idempotent. Then for any idempotent $e\in K_\lambda t(G)$ there
exists an infinite set of idempotents $e_1=e, e_2, \ldots $ such
that}
$$
                    e_ke_{k+1}=e_{k+1} \qquad (k\in \Bbb N). \tag9
$$
\endproclaim
{\bf Proof.} Suppose $K_\lambda t(G)$ does not contain a
nontrivial minimal idempotent. First we prove that for any
idempotent there exists an infinite set of idempotents $e_1, e_2,
\ldots$ in $K_\lambda t(G)$ satisfying condition (9).

Let $e_1$ be an idempotent of $K_\lambda t(G)$ and
$H_1=\gp{supp(e_1)}$. Then the ideal $K_\lambda t(G)e_1$ is not
minimal and contains a nontrivial ideal $\frak I_1$ of $K_\lambda
t(G)$. Let $0\not=x_1\in \frak I_1$ and $H_2=\gp{H_1, supp(x_1)}$.
Then $\overline{\frak I_1}=\frak I_1\cap K_\lambda H_2$ is an
ideal of $K_\lambda H_2$ and $\overline{\frak I_1}$ is generated
by the idempotent $e_2$ because $H_2$ is a finite subgroup of
$t(G)$ and the commutative algebra $K_\lambda H_2$ is semiprime.
It is easy to see that $e_1=e_2+f$, $f\not=0$ and $e_1e_2=e_2$.
Indeed, if $f=0$, then $e_1=e_2$ and $K_\lambda t(G)e_1=K_\lambda
t(G)e_2\subset\frak I_1$, which does not hold.  The ideal
$K_\lambda t(G)e_2$ contains a nontrivial ideal $\frak I_2$ of
$K_\lambda t(G)$. We choose a nonzero element $0\not=x_2\in \frak
I_2$ and construct the subgroup $H_3=\gp{H_2, supp(x_2)}$. The
ideal $\overline{\frak I_2}=\frak I_2\cap K_\lambda H_3$ is
generated by the idempotent $e_3$ and $e_2e_3=e_3\not= e_2$. This
method enables us to construct an infinite number of idempotents
$e_1, e_2, \ldots,  $ satisfying condition (9), which completes
the proof. \hfill \qed

\proclaim{Lemma 9} {\it Let $K$ be a field such that $char(K)$ does
not divide the order of any element of $t(G)$, and $U(K_\lambda
G)$ an $FC$-group. If the commutative algebra $K_\lambda t(G)$
contains an infinite number of idempotents $f_1, f_2, \ldots$ and
$g=[a, b]$ \quad  $(a, b\in G)$ is an element of order $n$ then
the commutators $[u_a,u_b]$ and $[a, b]$ have the same order and
$$
(f_i-f_j)(1-[u_a, u_b])=0\tag10
$$
for some $i\not=j$. }
\endproclaim
{\bf Proof.}
Let $g=[a,b]\not=1$, where  $a, b\in G$. By B.H.
Neumann's Theorem $G/t(G)$ is abelian, thus $g\in t(G)$ and
$1-f_i+u_bf_i$ is a unit in $K_\lambda G$. The element $u_a$ has a
finite number of conjugates in $U(K_\lambda G)$ and
$$
(1-f_i+u_b^{-1}f_i)u_a(1-f_i+u_bf_i)=u_a(1-f_i+[u_a, u_b]f_i).
$$
Consequently there exist $i$ and $j$ ($i<j$), such that
$$
1-f_i+[u_a, u_b]f_i=1-f_j+[u_a, u_b]f_j
$$
and
$$
(f_i-f_j)(1-[u_a,u_b])=0. \tag11
$$
If $n$ is the order of $g=[a, b]$ then
$$
[u_a,u_b]=\lambda_{a^{-1}, a}^{-1}\lambda_{b^{-1}, b}^{-1}
\lambda_{a^{-1}, b^{-1}}\lambda_{a^{-1}b^{-1}, a}
\lambda_{a^{-1}b^{-1}a, b}u_g
$$
and $[u_a, u_b]^n=\gamma\in U(K)$.
Then by (11) we have that $\gamma (f_i-f_j)=f_i-f_j$. Thus $\gamma=1$
and $[u_a, u_b]^n=1$.
\hfill \qed

\proclaim{Theorem 4}
{\it Let $K_\lambda G$ be an infinite twisted group algebra, and
$char(K)$ does not divide the order of any element of $t(G)$. If
$K_\lambda t(G)$ contains only a finite number of idempotents then
$U(K_\lambda G)$ is an $FC$-group if and only if $G$ is an $FC$-group
and the following conditions are satisfied:
\item{1.} all idempotents of $K_\lambda t(G)$ are central in $K_\lambda G$;
\item{2.} $\{\lambda_{h, h^{-1}}^{-1} \lambda_{h^{-1}, g}\lambda_{h^{-1}g, h}
\mid h\in H \}$ is a finite set for every $g\in G$;
\item{3.} $K_\lambda t(G)$ is a direct sum of a finite number of fields;
\item{4.} if $K_\lambda t(G)$ is infinite then it is central in $K_\lambda G$.}
\endproclaim
{\bf Proof. Necessity.}
By Lemma 4, 6 and 7 $K_\lambda t(G)$ is commutative, $G$ is an
$FC$-group and all idempotents of $K_\lambda t(G)$ are central in
$K_\lambda G$. Since $\{u_g^{-1}u_hu_g \mid g\in G\}$ is a finite
set, condition 2. of the Theorem is satisfied.

Since $K_\lambda t(G)$ contains only a finite number of idempotents (by Lemma
7) $K_\lambda t(G)$ is a direct sum of a finite number of fields. Let
$K_\lambda t(G)$ be infinite and $K_\lambda t(G)e_i$ is invariant under the
inner automorphism $\varphi(x)=u_g\m1xu_g$ for any $g\in G$.
Since $\gp{u_g, K_\lambda t(G)e_i\setminus \{0\} }$ is an $FC$-group  there are
exists an infinite subfield $L_g$ of $K_\lambda t(G)e_i$ such that $yu_g=u_gy$
for every $y\in L$.  Let $H=\gp{g, t(G)}$. Then $K_\lambda H$ is a subalgebra of
$K_\lambda G$ and (by Lemma 6) $K_\lambda t(G)$ is central in $K_\lambda
H$.

{\bf Sufficiency.}
Let $K_\lambda t(G)$ be a direct sum of fields,
$$
                    K_\lambda t(G)=F_1\oplus F_2\oplus \cdots \oplus F_t.
$$
Then $F_i=K_\lambda t(G)e_i$, where $e_i$ is a central idempotent
in $K_\lambda G$. It is easy to see that $K_\lambda G$ is a direct
sum of ideals
$$
     K_\lambda G=K_\lambda Ge_1\oplus \cdots \oplus K_\lambda Ge_t.       \tag12
$$
Let us prove that $K_\lambda Ge_q$ is isomorphic to a crossed product $F_q*H$
of the group $H=G/t(G)$ and the field $F_q$.

Let $R_l(G/t(G))$ be a fixed set of representatives of all left
cosets of the subgroup $t(G)$ in $G$. Any element $x\in K_\lambda
Ge_q$ can be written as
$$
x=e_qu_{c_1}\gamma_1+\cdots+e_qu_{c_s}\gamma_s,
$$
where $\gamma_k\in K_\lambda t(G), c_k\in
R_l(G/t(G))$. If $c_ic_j=c_kh$ \quad $(h\in t(G))$\quad  then
$$
     u_{c_i}u_{c_j}=u_{c_ic_j}\lambda_{c_i, c_j}=u_{c_kh} \lambda_{c_i,
          c_j}=u_{c_k}u_h\lambda_{c_k, h}^{-1} \lambda_{c_i, c_j}.
$$
We will construct the crossed product $F_q*H$, where
$$
H=\{h_i=c_it(G) \mid c_i\in R_l(G/t(G))\}.
$$
Let $\alpha\in F_q$ and  $\sigma$ be a map from $H$ to
the group of automorphism $Aut(F_q)$ of the field $F_q$ such that $\sigma
(h_i)(\alpha)=u_{c_i}^{-1}\alpha u_{c_i}$ and let $\mu_{h_i,
h_j}=u_h\lambda_{c_k, h}^{-1} \lambda_{c_i, c_j}$.

Clearly, the set $\mu =\{\mu_{a, b}\in U(F_q) \mid a, b\in H\}$ of
nonzero elements of field $F_q$, satisfies
$$
         \mu_{a, bc}\mu_{b, c}=\mu_{ab, c} \mu_{a, b}^{\sigma (c)},
$$
and
$$
        \alpha^{\sigma (a)\sigma (b)}=\mu_{a, b}^{-1}
                  \alpha^{\sigma (ab)} \mu_{a, b},
$$
where $\alpha \in F_q$ and $a, b, c\in H$.

Then $F_q*H=\{\sum_{h\in H}w_h\alpha_h \mid\alpha_h\in F_q\}$ is a crossed
product of the group $H$ and the field $F_q$ and we have $w_{d_i}w_{d_j}=
w_{d_k}\mu_{d_i, d_j}$ and $\alpha w_{d_i}=w_{d_i}\alpha ^{\sigma (d_i)}$.

Clearly, $F_q*H$ and $K_\lambda Ge_q$ are isomorphic because
$$
u_{c_i}\alpha u_{c_j}=u_{c_i}u_{c_j} (u_{c_j}^{-1}\alpha
u_{c_j})=u_{c_k}\mu_{c_i, c_j} \alpha^{\sigma (c_j)}.
$$
We know \cite{5} that the group of units of the crossed product
$K*H$ of torsionfree abelian group $H$ and  the field $K$ consists
of the elements  $w_h\alpha $, where  $\alpha \in
U(K)$, $h\in H$. By (12) for every $y\in U(K_\lambda G)$,
$$
                 y=u_{c_1}\gamma_1+\cdots +u_{c_t}\gamma_t \quad \text{and}\quad
                         y^{-1}=u_{c_1}^{-1}\gamma_1^{-1} +\cdots +
                               u_{c_t}^{-1}\gamma_t^{-1},
$$
where $\gamma_1, \ldots,  \gamma_t$ are orthogonal elements.

Let $x=\delta_1u_{d_1}+\cdots +\delta_tu_{d_t}\in U(K_\lambda G)$. Then
$$
yxy^{-1}=u_{c_1}\gamma_1\delta_1u_{d_1}u_{c_1}^{-1}\gamma_1^{-1}+\cdots
+u_{c_t} \gamma_1 \delta_tu_{d_t}u_{c_t}^{-1} \gamma_t^{-1}.
$$
If $K_\lambda t(G)$ is infinite then $K_\lambda t(G)\subseteq
\zeta (K_\lambda G)$ and
$$
        yxy^{-1}=\sum_{i=1}^t\delta_iu_{c_i}u_{d_i}u_{c_i}^{-1}=
   \sum_{i=1}^t \delta_i\lambda_{c_i, c_i^{-1}}^{-1}\lambda_{c_i,d_i}
           \lambda_{c_id_i, c_i^{-1}} u_{c_id_i{c_i}^{-1}}.
$$
Since $G$ is an $FC$-group, by condition 2. of the Theorem,  $x$ has a
finite number of conjugates, so $U(K_\lambda G)$ is an $FC$-group.

If $K_\lambda t(G)$ is finite then $F_q$ is a finite field and
$$
\split y^{-1}xy&=\sum_{i=1}^t
\gamma_i^{-1}u_{c_i}^{-1}\delta_iu_{d_i}u_{c_i}\gamma_i\\
&=\sum_{i=1}^t \lambda_{c_i, c_i^{-1}}^{-1}\lambda_{c_i^{-1}, d_i}
\lambda_{c_i^{-1} d_i, c_i} \gamma_1^{-1}\delta_i^{\sigma
(c_i^{-1})} \gamma_i^{\sigma (c_i^{-1}d_i^{-1}c_i)
}u_{c_i^{-1}d_ic_i}. \endsplit
$$
Since  $G$ is an $FC$-group and $F_q$ is a finite field,  $x$ has a
finite number of conjugates, so $U(K_\lambda G)$ is an $FC$-group.
\hfill \qed

\proclaim{Theorem 5}
{ \it Let $K_\lambda G$ be an infinite algebra over field
$K$, and $char(K)$ does not divide the order of any element of $t(G)$. If the
algebra $K_\lambda t(G)$ contains an infinite number of idempotents then
$U(K_\lambda G)$ is an $FC$-group if and only if $G$ is an $FC$-group and the
following conditions are satisfies
\item{1.} $K_\lambda t(G)$ is central in $K_\lambda G$ and contains a minimal
idempotent;
\item{2.} $\{\lambda_{h, h^{-1}}^{-1} \lambda_{h^{-1}, g}\lambda_{h^{-1}g, h}
\mid h\in H \}$ is a finite set for any $g\in G$;
\item{3.} the commutator subgroups of $G$ and of $\overline{G}$ are
isomorphic and $G'$ is either a finite group or isomorphic to the
group $\Bbb Z(q^\infty)$ $(q\not= p)$, and there exists an $n\in \Bbb
N$, such that the field $K$ does not contain the primitive $q^n$-th
root of $1$;
\item{4.} for every finite subgroup $H$ of the commutator subgroup of
$\overline{G}$ the element $e_H=\frac1{\mid H\mid }\sum_{h \in H}h$
is a nonzero idempotent of $K_\lambda t(G)$, and $K_\lambda
t(G)(1-e_H)$ is a direct sum of a finite number of fields.
\footnote{ If $K_\lambda G$ is a group ring, then 1) and 3) implies 4)
(see \cite{6} p.690, Lemma 4.3, also \cite{10}).}

}
\endproclaim
{\bf Proof. Necessity.}
By Lemma 4, 6 and 7 $K_\lambda t(G)$ is commutative, $G$ is an
$FC$-group and all idempotents of $K_\lambda t(G)$ are central in
$K_\lambda G$.

Let us prove that $K_\lambda t(G)$ contains a minimal idempotent.
Suppose the contrary. Let $a, b\in G$ and $1\not= [a, b]=g$. Since
$g$ is an element of finite order $n$, by Lemma 9,  $[u_a, u_b]^n=1$
and
$$
             f=\textstyle\frac1{n}(1+[u_a, u_b]^1+[u_a, u_b]^2+\cdots +[u_a, u_b]^{n-1})
$$
is an idempotent. By Lemma 11, for the idempotent $1-f$ one can
construct an infinite sequence of idempotents $e_1=1-f, e_2,
\ldots$\quad  satisfying (9). By Lemma 9,
$$
(1-[u_a, u_b])(e_i-e_j)=0,
$$
where $i<j$. Consequently $([u_a, u_b])^k(e_i-e_j)=(e_i-e_j)$ for
all $k$ and $f(e_i-e_j)=(e_i-e_j)$. This implies
$(1-f)(e_i-e_j)=0$.  Since $e_1=1-f$, $e_1(e_i-e_j)=0$. If we
multiply this equality from the right by the elements $e_2,
\ldots, e_{j-1}$, by (9) we obtain $(e_{j-1}-e_j)=0$. Now we
arrived at a contradiction, which proves that $K_\lambda t(G)$
contains a minimal idempotent.

It is easy to see that $t(G)$ is infinite, otherwise $K_\lambda t(G)$ would
contain a finite number of idempotents. $Kt(G)$ contains a minimal
idempotent $e$, and there exists only a finite number of elements $g\in t(G)$,
such that $eu_g=e$. Consequently $K_\lambda t(G)e$ is an infinite field
and contains $K$ as a subfield. Then as in the proof of Theorem 4, $K_\lambda t(G)$ is
central in $K_\lambda G$.

Since $\{ u_g^{-1}u_hu_g \mid g\in G \}$ is a finite set, we obtain condition
2. of the Theorem.

Suppose $c\in G'$ and
$$
                   c=[a_1, b_1][a_2, b_2]\cdots [a_n, b_n].
$$
Since $K_\lambda t(G)$ is central in $K_\lambda G$ and $1-e_i+e_iu_{b_k}\in
U(K_\lambda t(G))$ we have
$$
\split
\prod_{k=1}^n (1-e_i+e_iu_{b_k}^{-1}) u_{a_k}(1-e_i+e_iu_{b_k}) &=
\prod_{k=1}^n (u_{a_k}(1-e_i+e_i[u_{a_k}, u_{b_k}]))\\
&=\prod_{k=1}^n(u_{a_k})( \prod_{ k=1}^n (1-e_i+e_i[u_{a_k}, u_{b_k}])).
\endsplit
$$
for all $i\in \Bbb N$. Since each $u_{a_1}, u_{a_2}, \ldots,
u_{a_k}$ has a finite number of conjugates, there are only a
finite number of different elements of form $\prod_{ k=1}^n
(1-e_i+e_i[u_{a_k}, u_{b_k}])$. We denoted these elements by $w_1,
\ldots,  w_t$. Let
$$
W_r(c)=\{ i\in \Bbb N \mid \prod_{k=1}^n (1-e_i+e_i[u_{a_k}, u_{b_k}])=w_r \}.
$$
It is easy to see that the set of natural numbers $\Bbb N$ can be
written as a union of subsets $W_i(c)$ $( i=1, \ldots, t )$, of which at
least one is infinite. If $W_1(c)$ is infinite and $i, j\in W_1(c)$
then
$$
              (e_i-e_j)(1-\prod_{k=1}^n [u_{a_k}, u_{b_k}])=0. \tag13
$$
It implies that if
$$
\prod_{k=1}^n [u_{a_k}, u_{b_k}]=\gamma\in U(K),
$$
then $\gamma =1$.

Now we prove that the commutator subgroups of $G$ and of
$\overline{G}$ are isomorphic. It is easy to see that the map
$\tau (\lambda u_g)=g$ \quad  $(\lambda \in U(K), g\in G)$\quad
is a homomorphism from $\overline{G}$ to $G$. Every element $h\in
\overline{G}'$ can be written as
$$
            h=[u_{a_1}, u_{b_1}][u_{a_2}, u_{b_2}]\cdots [u_{a_n}, u_{b_n}].
$$
As we have shown above, if $h=\lambda \in U(K)$ then $\lambda=1$. Thus,
$\tau$ is an isomorphism from $\overline{G}'$ to $G'$.

Let $H$ be a finite subgroup of $\overline{G}'$. Then
$e_H=\frac1{\mid H\mid} \sum_{h\in H}h$ is an idempotent of
$K_\lambda t(G)$. Suppose that $K_\lambda t(G)(1-e_H)$ contains
infinite number of idempotents $e_1, e_2, \ldots  $. If $H=\{h_1,
h_2, \ldots,  h_s\}$, then as it is shown above, for every $h_j\in
H$,
$$
                     \Bbb N=W_1(h_j)\cup \cdots  \cup W_{r(j)}(h_j),
$$
where $j=1, 2, \ldots,  s$, and for every $k\not= l$, $W_k(h_j)$
and $W_l(h_j)$ have the empty intersection.

It is clear that there exists an infinite subset $M=W_{i_1}(h_1)\cap
\cdots \cap W_{i_s}(h_s)$. If $i, j\in M$, then by (13), we have
$(e_i-e_j)(1-h_r)=0$ for any $r$. Then
$$
e_i-e_j=\frac1{\mid H\mid}\sum_{r=1}^s (h_r(e_i-e_j))=e_H(e_i-e_j). \tag14
$$
Since $(e_i-e_j)\in K_\lambda t(G)(1-e_H)$, by (14),
$$
               (e_i-e_j)=(1-e_H)(e_i-e_j)= (e_i-e_j)-e_H(e_i-e_j)=0.
$$
Thus, $K_\lambda t(G)(1-e_H)$ contains a finite number of idempotents,
and by Lemma 7, it can be given as a direct sum of a finite number of fields.

Let us prove that there exist only finitely many elements of prime order in
$G'$.

Suppose the contrary. If $a, b\in G$ then $1\not= [a, b]=g\in
t(G)$. As we have seen above, if $h\in G'$, then there exists $\mu
\in U(K)$ such that the order of the element $\mu u_h$ equals the
order of $h$. Then there exists a countably infinite subgroup $S$,
generated by elements of prime order, such that $\gp{g}\cap S=1$.
By Pr\"uffer's Theorem \cite{9} $S$ is a direct product of cyclic
subgroups $S=\prod_{j}\gp{a_j}$ and $q_j$ is the order of element
$a_j$. Then
$$
e_j=\textstyle\frac1{q_j}(1+\mu u_{a_j}+(\mu u_{a_j})^2+\cdots
+(\mu u_{a_j})^{q_j-1})
$$
is a central idempotent and $x_i=1-e_i+e_iu_a\in U(K_\lambda G)$.
By lemma 9 $(e_i-e_j)(1-\mu u_g)=0$. Since $g\notin S$ we have  $i=j$, which does not hold.
Consequently $G'$ contains
only a finite number of elements of prime order and satisfies the
minimum condition for subgroups (see \cite{8}). Then
$$
G'\cong P_1\times P_2\times\cdots \times P_t\times H,
$$
where $P_i=\Bbb Z(q^\infty)$ and $\mid H\mid <\infty$. Let us prove
that either $G'=\Bbb Z(q^\infty)$ or $\mid G'\mid $ is finite.

Let $a, b\in G$ and $1\not= [a, b]=g\in t(G)$. Suppose there
exists $i$ such that $g\notin P_i= \gp{a_1, a_2, \ldots \mid
a_1^q=1, a_{j+1}^q=a_j}$. Then
$$
e_k=\textstyle\frac1{q^k}(1+\mu u_{a_j}+(\mu u_{a_j})^2+\cdots
+(\mu u_{a_j})^{q^k-1})
$$
is an idempotent and  $(e_i-e_j)(1-\mu u_g=0$. This is true only for
$i=j$, if $g\notin P_i$, which is
impossible. Thus, either  $G'\cong \Bbb Z(q^\infty)$ or $G'$ is a finite
subgroup.

Let $K$ be a field which contains primitive $q^n$-th root
$\varepsilon_n$ of $1$ for all $n$ and
$$
      P_1=\gp{a_1, a_2, \ldots  \mid a_1^{q^n}=1,\quad  a_{j+1}^{q^n}=a_j}.
$$
Put
$$
e_j=\textstyle\frac1{q^j}(1+\varepsilon_j\mu
u_{a_j}+(\varepsilon_j\mu u_{a_j})^2+\cdots +(\varepsilon_j\mu
u_{a_j})^{q^j-1}).
$$
If $i\not=j$ then  $(e_i-e_j)(1-\mu u_g)\not= 0$ and by
Lemma 9 this is impossible. Thus there exists $n\in \Bbb N$ such
that $K$ does not contain a primitive $q^n$-th root
$\varepsilon_n$ of $1$.

\smallskip
{\bf Sufficiency. } Let us prove that any element $u_g$ ($g\in G$) has a
finite number of conjugates in $U(K_\lambda G)$.

Let $\overline{G}=\{\kappa u_a\mid \kappa
\in U(K), a\in G\}$. We prove that $H=\gp{[u_g, \overline{G}]}$ is a finite subgroup in
$\overline{G}'$. If $\overline{G}'$ is finite, it is then obvious. If
$\overline{G}'$ is infinite then it is isomorphic to a subgroup of the group
$\Bbb Z(q^\infty)$. Any element of the group $\overline{G}$ is of the form
$\mu u_h$ \quad  $(\mu\in U(K), h\in G)$\quad  and
$$
             [u_g, \mu u_h]=\lambda_{g, g^{-1}}^{-1}\lambda_{h,
          h^{-1}}^{-1}\lambda_{g^{-1}, h^{-1}} \lambda_{g^{-1}h^{-1},
              g}\lambda_{g^{-1}h^{-1}g, h} u_{g^{-1}h^{-1}gh}.
$$
Since $G$ is an $FC$-group, and for a fixed element $g$ the set
$\{\lambda_{h, h^{-1}}^{-1} \lambda_{h^{-1}, g}\lambda_{h^{-1}g,
h} \mid h\in H \}$ is finite, the number of commutators $[u_g, \mu
u_h]$ is finite. These commutators generate a finite cyclic
subgroup $H$ of $\Bbb Z(q^\infty)$. The element $e_H=\frac1{\mid
H\mid} \sum_{h\in H} h$ is a nonzero idempotent in $K_\lambda
t(G)$ and by condition 4.  of the Theorem $K_\lambda t(G)(1-e_H)$
is a direct sum of a finite number of fields $K_\lambda
t(G)(1-e_H)f_i$ \quad  $(i=1, \ldots,  s)$.

In $K_\lambda t(G)$ we have the decomposition
$$
K_\lambda t(G)=K_\lambda t(G)e_H\oplus K_\lambda t(G)f_1\oplus \cdots \oplus
K_\lambda t(G)f_t.
$$
Then
$$
K_\lambda G=K_\lambda Ge_H\oplus K_\lambda Gf_1\oplus \cdots \oplus K_\lambda Gf_t.
$$
If $x\in U(K_\lambda G)$ then
$$
x=xe_H+xf_1+\cdots +xf_t\quad \text{and}\quad    x^{-1}=x^{-1}e_H+x^{-1}f_1+\cdots +x^{-1}f_t.
$$
Consequently
$$
        x^{-1}u_gx=x^{-1}e_Hu_gxe_h +x^{-1}f_1u_gxf_1 +\cdots +x^{-1}f_tu_gxf_t.
$$
We show that the element $xe_H$ is central in $U(K_\lambda G)$. If
$x=\alpha_1u_{h_1}+\cdots +\alpha_tu_{h_t}$, then
$$
\split u_gxe_H&=\alpha_1u_gu_{h_1}e_H+\cdots
+\alpha_tu_gu_{h_t}e_H\\
&=\alpha _1u_{h_1}u_g[u_g, u_{h_1}]e_H+\cdots
+\alpha_tu_{h_t}[u_g, u_{h_t}]e_H \endsplit
$$
and $[u_g, u_h]\in H$. Clearly, $[u_g, u_{h_k}]e_H=e_H$ and
$$
u_gxe_H=\alpha_1 u_{h_1}u_ge_H+\cdots  +\alpha_t u_{h_t}e_H= xe_Hu_g.
$$
$K_\lambda Gf_i$ is a crossed product $F*H$ of the group $H=G/t(G)$ and the
field $F=K_\lambda t(G)f_i$. We know (\cite{5}) that the group of units of
the crossed product $F*H$ of a torsionfree abelian group $H$ and a field $F$
consists of the elements $\alpha u_h$ \quad   $(\alpha \in U(F), h\in H)$.
The unit element $xf_i$ can be given as $\alpha_iu_{h_i}$,  where $h\in
G$ and  $\alpha_i$ is central in $U(K_\lambda Gf_i)$. Thus
$$
x^{-1}f_iu_gxf_i=\-u_{h_i}^{-1}\alpha_i^{-1}u_g\alpha_iu_{h_i}=\-
u_{h_i}^{-1}u_gu_{h_i}=\- \lambda_{h_i^{-1}, h_i}^{-1}
\lambda_{h_i^{-1}, g}\lambda_{h_i^{-1}g, h_i} u_{h_i^{-1}gh_i}.
$$
Therefore
$$
x^{-1}u_gx =u_g + \sum_{i=1}^t \lambda_{h_i^{-1}, h_i}^{-1}
\lambda_{h_i^{-1}, g} \lambda_{h_i^{-1}g, h_i} u_{h_i^{-1}gh_i}.
$$
Since $G$ is an $FC$-group, by condition 2.  of the Theorem, $u_g$
has a finite number of conjugates in $U(K_\lambda G)$.
\hfill \qed


\Refs

\ref\no 1
\by S.K.Sehgal, H.J. Zassenhaus
\paper Group rings whose units form an $FC$-group
\jour Math. Z.
\vol 153
\yr 1977
\pages 29-35
\endref

\ref\no 2
\by S.K.Sehgal, H.Cliff
\paper Group rings whose units form an $FC$-group
\jour Math. Z.
\vol 161
\yr 1978
\pages 169-183
\endref


\ref\no 3
\by B.H.Neumann
\paper Groups with finite classes of conjugate elements
\jour Proc. London Math.Soc.
\vol 1
\yr 1951
\pages 178-187
\endref

\ref\no 4
\by A.Kert\'esz
\book Lectures on artinian rings
\publ Akad\'emiai kiad\'o
\publaddr Budapest
\yr 1987
\endref

\ref\no 5
\by A.A.Bovdi
\book Group rings
\publ Kiev
\publaddr UMK VO
\yr 1988
\endref

\ref\no 6
\by D.S.Passman
\book The algebraic structure of group rings
\publ A Wiley-interscience publ
\publaddr New York- Syd\-ney- To\-ron\-to
\yr 1977
\endref

\ref\no 7
\by W.R.Scott
\paper On the multiplicative group of a division ring
\jour Proc. Amer. Soc
\vol 8
\yr 1957
\pages 303-305
\endref

\ref\no 8
\by L.Fuchs
\book Abelian groups
\publ Budapest:
\publaddr Publishing house of Hung. Acad. Sci.
\yr 1959
\endref

\ref\no 9
\by A.G.Kurosh
\book Theory of Groups
\publ New York:
\publaddr Chelsea
\yr 1955
\endref

\ref\no 10
\by J.S.Richardson
\paper Primitive idempotents and the socle in group rings of periodic
  abelian groups
\jour Composito mathematica Fasc
\vol 32
\yr 1976
\pages 203-223
\endref

\endRefs
\enddocument